\documentclass{amsart}
\usepackage{graphicx}
\usepackage{amsrefs,mathrsfs,math}
\newcommand\K{{\mathbb K}}
\numberwithin{equation}{section}

\begin{document}
\title{On Amenability of Group Algebras, II: graded algebras}
\author{Laurent Bartholdi}
\date{typeset \today; last timestamp 20070902}
\address{\'Ecole Polytechnique F\'ed\'erale de Lausanne (EPFL),
  Institut de Math\'ematiques B (IMB), 1015 Lausanne, Switzerland}
\email{laurent.bartholdi@gmail.com}
\begin{abstract}
  We show that, in a finitely generated amenable group $G$ with lower
  central series $(\gamma_n(G))$, the function
  $n\mapsto\rank(\gamma_n(G)/\gamma_{n+1}(G))$ grows subexponentially.

  This paper continues~\cite{bartholdi:aa1}'s study of amenability of
  affine algebras (based on the notion of almost-invariant
  finite-dimensional subspace), and applies it to graded algebras
  associated with finitely generated groups.

  We consider the graded deformation associated with the filtration of
  $\Bbbk G$ by powers of its augmentation ideal, and show that it has
  subexponential growth if $G$ is amenable. This yields the statement
  in the first paragraph, and answers a question by
  Vershik~\cite{vershik:amenability}, and another one by de la
  Harpe~\cite{harpe:uniform}.

  We also consider the graded deformation of a group ring $\Bbbk G$
  associated with a metric on $G$, and note that this deformation is
  amenable whenever $G$ is --- but also if $G$ has ``dead ends'' in
  its Cayley graph.
\end{abstract}
\maketitle

Warning! This paper contains a mistake in Lemma 4.3, and its main
claims should not be considered as proven!

\section{Introduction}
Throughout this paper, let $\Bbbk$ denote a commutative ring, and for
a $\Bbbk$-module $M$ let $\rank(M)$ denote its minimal number of
generators.  Recall from~\cite{bartholdi:aa1} the notion of amenable
algebra, which first appeared (in a slightly different form)
in~\cites{elek:amenaa,gromov:topinv1} and in~\cite{vaillant:folner} in
the context of $C^*$-algebras:

\begin{defn}\label{def:aa}
  Let $R$ be a $\Bbbk$-algebra and let $M$ be a right $R$-module. It
  is \emph{amenable} if, for every $\epsilon>0$ and every finite-rank
  subspace $S$ of $R$, there exists a finite-rank subspace $F$ of $M$
  such that
  \[\frac{\rank(F+FS)-\rank(F)}{\rank(F)}<\epsilon;\]
  any such $F$ is called \emph{$(S,\epsilon)$-invariant}. 

  $M$ is \emph{exhaustively amenable} if furthermore in the definition
  above the space $F$ may be required to contain any specified
  finite-rank subspace.
\end{defn}

This is a counterpart to amenability of $G$-sets for a group $G$,
where according to F\o lner's definition~\cite{folner:banach} the
$G$-set $X$ is \emph{amenable} if, for every $\epsilon>0$ and every
finite subset $S\subset G$, there exists a finite subset $F\subset X$
such that $(\#(F\cup FS)-\#F)/\#F<\epsilon$; such $F$ being also
called $(S,\epsilon)$-invariant.  The main result
of~\cite{bartholdi:aa1} was:
\begin{thm}\label{thm:main}
  Let $\K$ be a field, and let $X$ be a $G$-set. Then $X$ is amenable
  if and only if its linear envelope $\K X$ is amenable.
\end{thm}


\subsection{Growth of modules}\label{ss:growth}
Let $R$ be an \emph{affine} algebra, i.e.\ a finitely generated
associative algebra over a commutative ring $\Bbbk$, generated by the
finite-rank subspace $S$.  Let $M$ be a finitely-generated $R$-module,
generated (as an $R$-module) by a finite-rank subspace $T$. Let us
write here and below
\[TS^{\le n}=T+TS+TS^2+\dots+TS^n.\] The module $M$ is then filtered
by the exhausting sequence of subspaces $T\le TS^{\le1}\le
TS^{\le2}\le\cdots\le TS^{\le n}\le\cdots$, and the \emph{growth} of
$M$ is the sequence $(r_n)_{n\ge1}$ defined by $r_n=\rank(TS^{\le n})$. We
say that $M$ has \emph{subexponential growth} if
$\lim\sqrt[n]{r_n}=1$. This property does not depend on the choices of
$S$ and $T$.
\begin{prop}\label{prop:subam}
  If $M$ has subexponential growth, then $M$ is exhaustively amenable.
\end{prop}
\begin{proof}
  Given $U\le M$ of finite rank and $\epsilon>0$, let $d\in\N$ be such
  that $U\le TS^{\le d}$. Since $(r_n)$ grows subexponentially,
  $\liminf r_{n+1}/r_n=1$, so there exists $n>d$ such that
  $r_{n+1}/r_n<1+\epsilon$.  Set $F=TS^{\le n}$. Then $F+FS\le TS^{\le
    n+1}$, so $(\rank(FS)-\rank(F))/\rank(F)<\epsilon$ and $F$
  contains $U$.
\end{proof}

We consider in~\S\ref{ss:crystal} the graded module $M_0=\bigoplus_n
TS^{\le n}/TS^{\le n-1}$ associated with the ascending filtration of
$M$ by its finite-rank subspaces $TS^n$. If we consider $R=\Bbbk G$ a
group ring, $S$ the $\Bbbk$-span of a finite generating set of $G$,
and $M=R$ with $T=\{1\}$, then $M_0=R_0$ is an algebra with basis in
bijection with $G$, and with product derived from multiplication in
$G$ except that $g\cdot h=0$ if the length of $gh$ is strictly less
than the sum of the lengths of $g$ and $h$. We prove:
\begin{thm}\label{thm:dead}
  Let $G=\langle S\rangle$ be a finitely generated group.  If $G$ is
  amenable, or has dead ends, then $(\Bbbk G)_0$ is amenable.  If $G$
  is amenable, or has infinitely many dead ends, then $(\Bbbk G)_0$ is
  exhaustively amenable.
\end{thm}

The proof of Theorem~\ref{thm:dead} relies on fact that $(\Bbbk
G)_0$ has a monomial basis, namely a basis $B$ with $B\cdot B\subseteq
B\cup\{0\}$, this prompts the
\begin{question}
  If $M$ is an amenable finitely-generated $R$-module, is then $M_0$
  also amenable?
\end{question}

\subsection{Augmented algebras}\label{ss:augmented}
Assume now that $R$ has an augmentation\footnote{Namely, a morphism of
  unital algebras to the scalars $\Bbbk$.}$\varepsilon:R\to\Bbbk$,
with kernel $\varpi$; as a typical example $R=\Bbbk G$ and
$\varpi=\langle g-1:\,g\in G\rangle$. Let $M$ be an $R$-module. Then
$M$ admits a descending filtration $M\ge M\varpi\ge
M\varpi^2\ge\cdots\ge M\varpi^n\ge\cdots$; and an associated graded
module \[\overline M=\bigoplus_{n\ge 0}M\varpi^n/M\varpi^{n+1}.\] In
particular, $R$ admits an associated graded algebra $\overline
R=\bigoplus_{n\ge0}\varpi^n/\varpi^{n+1}$. If $M$ is generated (as an
$R$-module) by a subspace $T$, then $\overline M$ is a graded
$\overline R$-module, generated by the image $\overline T$ of $T$ in
$M/\varpi M$.

If $R$ is affine, say generated by a finite-rank subspace $S$, then
$\overline R$ is also affine: the projection $\overline
S=\{s-\varepsilon(s)+\varpi^2:\,s\in S\}$ of $S$ in $\varpi/\varpi^2$
actually equals $\varpi/\varpi^2$ irrespective of $S$, and generates
$\overline R$.  The \emph{growth} of $\overline M$ is then its growth
in the sense of~\S\ref{ss:growth}, with respect to the generating
subspace $\overline S$ of $\overline R$ and $\overline T$ of
$\overline M$.

The main result of this paper is the following theorem, proven
in~\S\ref{ss:vershik}:
\begin{thm}\label{thm:vershik}
  If $G$ is an amenable, finitely generated group, then
  $\overline{\Bbbk G}$ has subexponential growth (and therefore is
  exhaustively amenable by Proposition~\ref{prop:subam}).
\end{thm}

Vershik conjectured in~\cite{vershik:amenability}*{page 326} that if
$G$ is an amenable, finitely generated group and if $\varpi$ denote
the augmentation ideal in $\Z G$, then $\rank(\varpi^n/\varpi^{n+1})$
grows subexponentially. This follows from Theorem~\ref{thm:vershik}
with $\Bbbk=\Z$.

Theorem~\ref{thm:vershik} actually follows from the following
statement: if $R$ is an amenable augmented affine algebra with a basis
consisting of invertible elements, then $\overline R$ has
subexponential growth. This prompts the
\begin{question}
  Let $R$ be an amenable augmented affine algebra. Does $\overline R$
  necessarily have subexponential growth?

  Let more generally $M$ be an amenable $R$-module, where $R$ is any
  augmented affine algebra. Does $\overline M$ necessarily have
  subexponential growth?
\end{question}


\subsection{Golod-Shafarevich groups}
Let $\K$ be a field of characteristic $p>0$, and assume that $G$ is a
residually-$p$ group. This means that the series $(G_{n,p})_{n\ge1}$
of \emph{dimension} subgroups defined by $G_{1,p}=G$ and
$G_{n+1,p}=[G_{n,p},G](G_{\lceil n/p\rceil,p})^p$ for $n\ge1$,
satisfies $\bigcap G_{n,p}=\{1\}$.  The \emph{degree} of $g\in G$,
written $\deg_p(g)$, is the maximal $n\in\N\cup\{\infty\}$ such that
$g\in G_{n,p}$.

It is well known~\cites{jennings:gpring,lazard:uea} that $\mathcal
L:=\bigoplus_{n\ge1}(G_{n,p}/G_{n+1,p})\otimes_{\F}\K$ has the structure of
a restricted Lie algebra over $\K$; and that $\overline{\K G}$ is the
universal restricted enveloping algebra of $\mathcal L$. We will use
this fact to relate Theorem~\ref{thm:vershik} to purely
group-theoretical statements, e.g.\ Corollary~\ref{cor:lcs}.

Golod constructed in~\cite{golod:nil} for every prime $p$ a finitely
generated infinite torsion $p$-group. His method is quite flexible,
and was generalized as follows:
\begin{defn}
  A group $G$ is a \emph{Golod-Shafarevich group} for the prime $p$ if
  it admits a presentation $G=F/\langle \mathscr R\rangle^F$ in which
  $F$ is a free group of rank $d$ and $\mathscr R\subset F$ is a set
  of relators, such that for some $t\in(0,1)$ we have
  \[1-dt+\sum_{r\in\mathscr R}t^{\deg_pr}<0,\]
  where $\deg_p$ denotes degree with respect to the filtration
  $(F_{n,p})_{n\ge1}$ of $F$.
\end{defn}
\begin{prop}[\cite{herstein:rings}; \cite{koch:galois};
  \cite{huppert-b:fg2}*{\S VIII.12}; \cite{bartholdi-g:lie}*{}]
  All Golod-Shafarevich groups are infinite. If $\varpi$ denote the
  augmentation ideal in $\F G$, then $\dim_{\F}(\varpi^n/\varpi^{n+1})$
  grows exponentially (at rate at least $1/t$).
\end{prop}
On the other hand, there are torsion groups that are Golod-Shafarevich
(this solved Burnside's problem~\cite{burnside:question}).

The second part of the following result
answers~\cite{harpe:uniform}*{Open Problem 5.2}
and~\cite{harpe:mfap}*{Question 7}:
\begin{cor}
  Golod-Shafarevich groups are not amenable. In particular, there
  exist non-amenable residually-$p$ torsion groups.
\end{cor}
Ershov~\cite{ershov:goshat} has constructed Golod-Shafarevich groups
that have property (T) --- and no infinite (T) group can be amenable.
Any Golod-Shafarevich group admits a Golod-Shafarevich torsion
quotient, which will still have property (T). This
answers~\cite{harpe:uniform}*{Open Problem 5.2} by a different method.

\subsection{Lower central series}
Theorem~\ref{thm:vershik} has the following purely group-theoretic
consequence, also proven in~\S\ref{ss:vershik}:
\begin{cor}\label{cor:lcs}
  Let $G$ be an amenable, finitely-generated group, and let
  $(\gamma_n(G))_{n\ge1}$ denote its lower central series. Then the
  function $n\mapsto\rank(\gamma_n(G)/\gamma_{n+1}(G))$ grows
  subexponentially.
\end{cor}
Note that this function may grow arbitrarily close to an exponential
function. Indeed Petrogradsky showed
in~\cite{petrogradsky:polynilpotent} that if $G$ be the free
$k$-generated soluble group of solubility class $q\ge3$, then
\[\rank(\gamma_n(G)/\gamma_{n+1}(G))\cong\exp\left(\bigg(\frac{(k-1)\zeta(k)}{\log\log\cdots n}\bigg)^{1/k}n\right),\]
with $q-3$ iterated logarithms in the expression above.

The converse of Theorem~\ref{thm:vershik} does not hold (see
Remark~\ref{rem:vershik}). However, the next statement trivially
follows from Theorem~\ref{thm:vershik}, and raises the question after
it:
\begin{cor}
  Let $G$ be an amenable group.  Then for every finitely-generated
  subgroup $H$ of $G$ the function
  $n\mapsto\rank(\gamma_n(H)/\gamma_{n+1}(H))$ grows subexponentially.
\end{cor}
\begin{question}
  Does there exist a non-amenable
  residually-nilpotent\footnote{Namely, such that
    $\bigcap\gamma_n(G)=1$.} group $G$ such that for every
  finitely-generated subgroup $H$ of $G$ the function
  $n\mapsto\rank(\gamma_n(H)/\gamma_{n+1}(H))$ grows subexponentially?
\end{question}

\subsection{Acknowledgments}
The author is grateful to Mikhail Ershov, Fran\c cois Gu\'eritaud,
Pierre de la Harpe and Fabrice Krieger for generous feedback and/or
entertaining and stimulating discussions.

\section{Hecke and crystal algebras}\label{ss:crystal}
We prove Theorem~\ref{thm:dead} in this section, phrasing it in a
slightly more general manner.  Let $G$ be a group with fixed
generating set $S$. Denote by $\ell(g)$ the length of $g\in G$ in the
word metric:
\[\ell(g)=\min\{n:\, g=s_1\dots s_n,\,s_i\in S\}.\]
Choose $\lambda\in\Bbbk$, and define the ``Hecke
algebra''\footnote{The terminology comes from the classical Hecke
  algebra associated with the symmetric group.}  $(\Bbbk G)_\lambda$
as follows: it is isomorphic to $\Bbbk G$ as a $\Bbbk$-module; it has
a basis $\{\delta_g\}_{g\in G}$; and multiplication is given by
\[\delta_g\delta_h=\lambda^{\ell(g)+\ell(h)-\ell(gh)}\delta_{gh}.\]
Note that, although the notation does not make it explicit, $(\Bbbk
G)_\lambda$ depends on the choice of $S$.

If $\lambda$ is invertible, then $(\Bbbk G)_\lambda$ is isomorphic to
$\Bbbk G$, the isomorphism being given by $\delta_g\mapsto
\lambda^{\ell(g)}g$. The universal cases $\Bbbk=\Z[\lambda]$ and
$\Q[[\lambda]]$ should deserve particular consideration.

Quite on the contrary, $(\Bbbk G)_0$ is a graded algebra (with degree
function $\ell$), which we call the \emph{crystal}\footnote{The
  terminology comes from statistical mechanics, where the parameter
  $\lambda$ of the deformation of $\Bbbk G$ is interpreted as
  temperature.} of $\Bbbk G$.  Note that $(\Bbbk G)_0$ is the
associated graded algebra $\bigoplus\Bbbk S^{\le n}/\Bbbk S^{\le n-1}$
of the filtered algebra considered in \S\ref{ss:growth}.  The first
part of Theorem~\ref{thm:dead} generalizes as:
\begin{prop}\label{prop:crystal}
  Let $G$ be an amenable group with fixed generating set. Then $(\Bbbk
  G)_\lambda$ is amenable for all $\lambda\in\Bbbk$.
\end{prop}
\begin{proof}
  Let $\epsilon>0$ be given, and let\footnote{In this proof, we use
    $S$ for a subset of a group, and $S'$ for a subspace of an
    algebra.} $S'$ be a finite-rank subspace of $(\Bbbk G)_\lambda$.
  Let $S$ denote the support of $S'$, i.e.\ the set of those $g\in G$
  such that $\delta_g$ has a non-zero co\"efficient in some element of
  $S'$; it is a finite subset of $G$.  Since $G$ is amenable, there
  exists a finite subset $F$ of $G$ with
  \begin{equation}\label{eq:crystal:1}
    (\#(F\cup FS)-\#F)/\#F<\epsilon.
  \end{equation}
  Set $F'=\bigoplus_{g\in F}\Bbbk\delta_g$, a finite-rank subspace of
  $(\Bbbk G)_\lambda$. We have $\rank F'=\#F$ and $F'S'\le\Bbbk(FS)$,
  so $\rank(F'S')\le\#FS$, and $\rank(F'+F'S')\le\#(F\cup FS)$, whence
  \begin{equation}\label{eq:crystal:2}
    \frac{\rank(F'+F'S')-\rank(F')}{\rank(F')}<\epsilon,
  \end{equation}
  so $(\Bbbk G)_\lambda$ is amenable.
\end{proof}

It is however possible for $(\Bbbk G)_0$ to be amenable, yet for $G$
not to be amenable. The example in Proposition~\ref{prop:bogo} appears
in~\cite{bogopolski:bilipschitz} (with a small typographical mistake).

\subsection{Dead ends} Say $g\in G$ is a \emph{dead end} if
$\ell(gs)\le\ell(g)$ for all $s\in S$. Note that, although the
notation does not make it apparent, this property strongly depends on
the choice of $S$. Zoran \v Suni\'k has informed me that every group
admits a generating set for which the group contains a dead end. The
second part of Theorem~\ref{thm:dead} reads:

\begin{prop}
  If $G$ has a dead end (with respect to $S$), then $(\Bbbk G)_0$ is
  amenable.

  If $G$ has infinitely many dead ends, then $(\Bbbk G)_0$ is
  exhaustively amenable.
\end{prop}
\begin{proof}
  Let first $g\in G$ be a dead end. Set $F=\Bbbk\delta_g$. Then
  $F+Fs=F$ for any $s\in(\Bbbk G)_0$, so $(\Bbbk G)_0$ is amenable.

  Let now $g_1,g_2,\dots$ be an infinite set of dead ends in $G$.
  Given $\epsilon>0$ and $E\le(\Bbbk G)_0$ of rank $n$, consider
  $F=E+\Bbbk\delta_{g_1}+\dots+\Bbbk\delta_{g_{\lceil
      n/\epsilon\rceil}}$.  Then $F$ contains $E$ and
  $\rank(F+Fs)\le\rank F+\rank E$ and $\rank F\ge\rank E/\epsilon$.
\end{proof}

\begin{prop}[Bogopolski]\label{prop:bogo}
  For all $k\ge 3$, the ``triangle'' group
  \[G_k=T_{3,3,k}=\langle x,y\mid x^3,y^3,(xy)^k\rangle\] contains
  infinitely many dead ends.
\end{prop}
\begin{proof}
  This group is hyperbolic; it acts by isometries on hyperbolic space
  and preserves the semiregular tiling of $\mathbb H^2$ by triangles
  (with edges labeled $x$ and $y$) and $2k$-gons (with edges labeled
  periodically $x,y$).  The $1$-skeleton of this tiling is the Cayley
  graph of $G_k$.

  Assume first that $k$ is even. Consider for all
  $n\in\Z\setminus\{0\}$ the element $d_n\in G_k$ defined by
  \[d_{2n}=((xy)^{k/2}(yx)^{k/2})^n,\qquad
  d_{2n+1}=((xy)^{k/2}(yx)^{k/2})^n(xy)^{k/2}.\]

  Consider furthermore the automorphism $\phi$ of $G_k$ defined by
  $x\mapsto y^{-1}$ and $y\mapsto x^{-1}$. It is easy to check that
  $d_n$ is fixed by $\phi$.

  Geometrically, $d_n$ is the hyperbolic isometry whose axis cuts
  through a doubly-infinite sequence of vertex-abutting $2k$-gons in
  their middle, and translates by $n$ of them; $\phi$ is the
  reflection through this axis.

  Let $s_1\dots s_m$ be a word of minimal length representing $d_n$,
  with $s_i\in S=\{x,y,x^{-1},y^{-1}\}$; so $\ell(d_n)=m$. Then
  $\ell(d_ns)\le m$ for all $s\in S$:
  \begin{itemize}
  \item if $s=s_m^{-1}$ this is clear;
  \item if $s=s_m$ then $d_ns=s_1\dots s_{m-1}s_m^{-1}$ whence
    $\ell(d_ns)\le m$;
  \item otherwise, $d_n^\phi=d_n=s_1^\phi\dots s_m^\phi$ also of
    minimal length, and $s_m^\phi\in\{s,s^{-1}\}$, so the previous
    cases apply.
  \end{itemize}

  If $k$ is odd, consider for all $n\ne0$ the element
  $d_n=((xy)^{(k-1)/2}x)^n$, with a similar geometric interpretation
  as above. The same arguments apply.
\end{proof}
It is in fact not hard to see that these are the only dead ends in
$G_k$. Since for $k\ge4$ the $G_k$ are non-elementary hyperbolic
groups, they are certainly not amenable.

There unfortunately does not seem to be any natural condition to
impose on $(\Bbbk G)_0$ to ensure that $G$ be amenable.

\section{Tileable amenable groups}\label{ss:tag}
We prove in this section a result by Weiss~\cite{weiss:monotileable},
based on earlier work by Ornstein and
Weiss~\cite{ornstein-weiss:entropyiso}*{\S I.2}. We follow the sketch
of a proof by Gromov~\cite{gromov:topinv1}*{pages~336--337}, adapting
it so as to prepare the ground for a generalization to modules
in~\S\ref{ss:tae}. We will not use the results in this \S, but present
them as a warm-up for \S\ref{ss:tae}.

Let $X$ denote a $G$-set. We defined amenability with respect to the
\emph{outer envelope} $AK$ of a subset $A\subset X$ with respect to
$K\subset G$. It will be useful in this section to consider, again for
$A\subset X$ and $K\subset G$, the \emph{inverse envelope}
\[AK^* := \{x\in X:\,xK\cap A\neq\emptyset\}.\]
The easy properties
\[A(K\cup L)^*=AK^*\cup AL^*,\qquad (A\cup B)K^*=AK^*\cup BK^*\]
follow immediately from $AK^*=A\{k^{-1}:\,k\in K\}$.

\begin{thm}[Weiss]\label{thm:weiss}
  Let $G$ be an amenable group, let $K\subseteq G$ be a finite subset,
  let $\epsilon>0$ be given, and let $N_0\ge N_1\ge\dots$ be a nested
  sequence of finite-index normal subgroups of $G$ such that
  $\bigcap_{n\in\N} N_n=\{1\}$.

  Then for all $n\gg0$ there exists a $(K,\epsilon)$-invariant subset
  $T_n\subseteq G$ that is a transversal for $N_n$ in $G$.
\end{thm}

Let $\zeta\ge1$ and $\delta\in(0,1)$ be constants to be fixed later.
In this section and the next we shall use, for $\mu\ge\delta$, the
transformation $\Theta_\mu:\R_+^2\to\R_+^2$ given by
\[\Theta_\mu(\nu,\alpha)=\Big(\nu+\mu(1-\alpha),\alpha+\frac{\mu(1-\alpha)}{1-\delta}\zeta\Big).\]
The proof of Theorem~\ref{thm:weiss} relies on the following
\begin{lem}\label{lem:weiss}
  Let $\Omega$ be a finite $G$-set; let $B$ be a subset of $\Omega$
  and let $K,L$ be finite subsets of $G$; assume that for all
  $x\in\Omega$ the orbit map
  \begin{equation}\label{eq:weiss:i}
    K\to\Omega,\; k\mapsto xk\qquad\text{is injective}.
  \end{equation}
  Let $\alpha,\nu\in[0,1)$ and $\zeta\ge1$ be such that
  \begin{xalignat*}{2}
    \#(KL^*)&\le\zeta\#K, & \#(BK^*)&\le\alpha\#\Omega,\\
    \#B&=\nu\#\Omega, & \#(BL^*)&\le\alpha\#\Omega;
  \end{xalignat*}
  Let furthermore $\delta\in(0,1)$ be given.

  Then there exist $s\ge 1$, elements $x_1,\dots,x_s\in\Omega$, and
  $\mu\ge\delta$ such that, setting $B_0=B$, $B_i=B_{i-1}\cup x_iK$
  for $1\le i\le s$ and $\Theta_\mu(\nu,\alpha)=(\nu',\alpha')$, we
  have
  \begin{gather}
    \label{eq:weiss:1}
    \#(x_iK\cap B_{i-1}) \le\delta\#K\text{ for all }i\in\{1,\dots,s\};\\
    \label{eq:weiss:2} \#B_s =\nu'\#\Omega;\\
    \label{eq:weiss:3} \#(B_sL^*) \le\alpha'\#\Omega.
  \end{gather}
\end{lem}
Informally, the lemma says that, given a set $K$ with small
$L$-boundary, and a set $B\subset\Omega$ with small $K$- and
$L$-boundary, one can construct a quantifiably larger set $B_s$ with
small $L$-boundary.

\begin{proof}
  Let $x_1,\dots,x_s\in\Omega$ be a maximal-length sequence of
  elements such that~\eqref{eq:weiss:1} holds.  Then, by maximality of
  $(x_1,\dots,x_s)$, for all $x\in\Omega$ we have
  \begin{equation}\label{eq:weiss:p1}
    \#(B_s\cap xK)=\#\{(b,k):b\in B_s,k\in K,xk=b\}>\delta\#K.
  \end{equation}
  We deduce
  \begin{align*}
    \#B_s\#K&=\#\{(b,k,x):\,b\in B_s,k\in K,x\in\Omega,xk=b\}\\
    &=\sum_{x\in\Omega}\#\{(b,k):b\in B_s,k\in K,xk=b\}
    =\sum_{x\in\Omega}\#(B_s\cap xK)\\
    &=\sum_{x\in BK^*}\#\{(b,k):b\in B_s,k\in K,xk=b\}+\sum_{x\in
      \Omega\setminus BK^*}\#(B_s\cap xK)\\
    &=\sum_{b\in B_s}\sum_{k\in K}\#\{BK^*\cap
    \{bk^{-1}\}\}+\sum_{x\in\Omega\setminus BK^*}\#(B_s\cap xK)\\
    &\ge\sum_{b\in B}\sum_{k\in K}\#\{BK^*\cap
    \{bk^{-1}\}\}+\sum_{x\in\Omega\setminus BK^*}\delta\#K\text{ by }\eqref{eq:weiss:p1}\\
    &\ge\#B\#K+(1-\alpha)\#\Omega\delta\#K.
  \end{align*}
  If we divide by $\#K\#\Omega$ and set
  $\mu=(\#B_s/\#\Omega-\nu)/(1-\alpha)$, then we get $\mu\ge\delta$,
  and therefore $s\ge1$.  Note that~\eqref{eq:weiss:2} holds by the
  choice of $\mu$. Next,
  \[\#B_s\ge\#B+\sum_{i=1}^s\#(x_iK\setminus
  B_{i-1})\ge\#B+s(1-\delta)\#K\text{ by }(\ref{eq:weiss:i},\ref{eq:weiss:1})
  \]
  so $s\#K\le(\nu'-\nu)\#\Omega/(1-\delta)$; and
  \[\#(B_sL^*)\le\#(BL^*)+\sum_{i=1}^s\#(x_iKL^*)\le\alpha\#\Omega+s\zeta\#K,
  \]
  from which~\eqref{eq:weiss:3} follows.
\end{proof}

\begin{rem}\label{rem:tiling}
  We shall apply Lemma~\ref{lem:weiss} using the following strategy:
  we construct a finite but very long sequence $K_1,\dots,K_t$ of F\o
  lner sets, each with a very small boundary with respect to its
  predecessors and to $K$. We then find a finite quotient $\Omega$ of
  $G$ in which these F\o lner sets embed. Then, starting from $K_t$
  down to $K_1$, we apply $t$ times Lemma~\ref{lem:weiss} to cover
  most of $\Omega$ by images of translates of the $K_i$, taking each
  time as many copies as possible subject to them having extremely
  small overlaps (see condition~\eqref{eq:weiss:1}). There remains a
  small part of $\Omega$ that is not covered. We then lift these
  images back to $G$, and lift the small remainder arbitrarily. We
  have obtained a transversal, consisting mostly of pieces carved out
  of F\o lner sets by other F\o lner sets.  See Figure~\ref{fig:1} for
  an illustration.
  \begin{figure}
    \[\includegraphics[height=15cm]{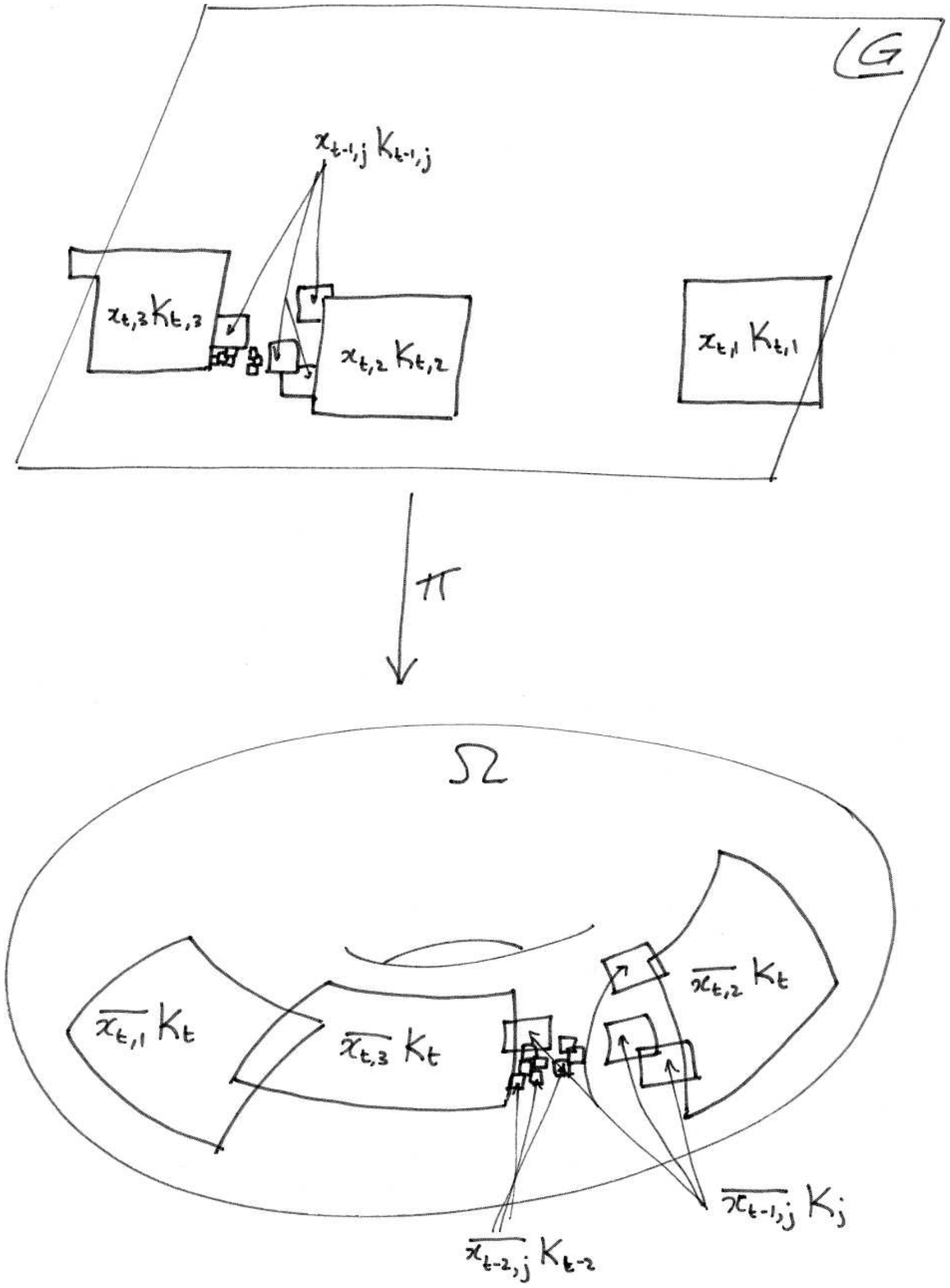}\]
    \caption{The proof of Theorem~\ref{thm:weiss}, illustrated in a picture}
    \label{fig:1}
  \end{figure}

  Note that both co\"ordinates of $\Theta_\mu$ increase monotonically
  in $\mu$. This property helps us formalize the intuitive notion that
  ``it cannot hurt us, in constructing the covering of $\Omega$, if
  the biggest F\o lner sets leave unexpectedly small holes between
  them, even though small holes are harder to cover efficiently with
  smaller F\o lner sets''.\hfill $\triangle$
\end{rem}

\begin{proof}[Proof of Theorem~\ref{thm:weiss}]
  We are given $K$ and $\epsilon$. Without loss of generality, we
  suppose $1\in K$. Fix in order the following data:
  \begin{enumerate}
  \item Choose $\delta>0$, the ``allowed overlap'' as
    in~\eqref{eq:weiss:1}, such that $\delta\#K<\frac\epsilon2$ and
    $(1+\frac\epsilon2)(1-\delta)>1$;
  \item Choose $\zeta>1$ such that
    $(1-\delta)/\zeta>1-\epsilon/(2\#K)$, the ``relative F\o lner
    constant''.
  \end{enumerate}
  For all $t\in\N$ set
  $(\overline{\nu_t},\overline{\alpha_t})=\Theta_\delta^t(0,0)$.  Then
  it is clear that
  $\overline{\nu_t}/\overline{\alpha_t}=(1-\delta)/\zeta$ for all
  $t>0$, and that $\limsup\overline{\alpha_t}\ge1$, so
  $\limsup\overline{\nu_t}\ge(1-\delta)/\zeta$.
  \begin{enumerate}\setcounter{enumi}{2}
  \item Let $t\in\N$, the ``tower height'', be such that
    $\overline{\nu_t}>1-\epsilon/(2\#K)$.
  \end{enumerate}
  Then for any $\mu_1,\dots,\mu_t\ge\delta$ we also have
  \[\Theta_{\mu_t}(\cdots\Theta_{\mu_1}(0,0)\cdots)\in(1-\epsilon/(2\#K),1]\times[0,1].\]
  \begin{enumerate}\setcounter{enumi}{3}
  \item Using amenability of $G$, construct finite subsets
    $K_0=K,K_1,\dots,K_t$ of $G$, the ``Rokhlin tower'', such that for
    all $j<i$ in $\{1,\dots,t\}$ one has $\#(K_iK_j^*)<\zeta\#K_i$,
    and such that
    \begin{equation}\label{eq:weiss:6}
      \#(K_iK)\le\big(1+\frac\epsilon2\big)(1-\delta)\#K_i \text{ for all
      }i\in\{1,\dots,t\}.
    \end{equation}
  \item Choose $n\in\N$, the ``quotient index'', large enough so that
    $K_iK_i^*\cap N_n=\{1\}$ for all $i\in\{1,\dots,t\}$.
  \end{enumerate}
  Set $\Omega=G/N_n$, let $\pi:G\to\Omega$ denote the natural quotient
  map, and remark, by the choice of $n$, that $K_i\to\Omega$,
  $k\mapsto xk$ is injective for all $i\in\{1,\dots,t\}$ and all
  $x\in\Omega$.
  
  Now start with $A_{t,0}=\emptyset$ and
  $(\nu_t,\alpha_t)=(0,0)=(\overline{\nu_0},\overline{\alpha_0})$, and
  apply $t$ times Lemma~\ref{lem:weiss}: at step $i=t,t-1,\dots,1$,
    apply it with $B:=A_{i,0}$, $K:=\pi(K_i)$, $L:=\pi(K_{i-1})$, and
    $(\nu,\alpha):=(\nu_i,\alpha_i)$. Lemma~\ref{lem:weiss} produces for
    us
    \begin{itemize}
    \item constants $s\ge1$ and $\mu\ge\delta$, which we rechristen
      $s(i)$ and $\mu_i$;
    \item a sequence $x_1,\dots,x_{s(i)}$ in $\Omega$, which we rechristen
      $\overline{x_{i,1}},\dots,\overline{x_{i,s(i)}}$;
    \item a sequence $B_0\subseteq\dots\subseteq B_{s(i)}$ of subsets of
      $\Omega$, which we rechristen
      $A_{i,0},\dots,A_{i,s(i)}=:A_{i-1,0}$.
    \end{itemize}
    Set $(\nu_{i-1},\alpha_{i-1})=\Theta_{\mu_i}(\nu_i,\alpha_i)$, and
    note that $\nu_{i-1}\ge\overline{\nu_{t-(i-1)}}$ by induction on $i$
    and because $\mu_i\ge\delta$ and $\Theta_\mu$ is monotonically
    increasing in $\mu$.  Finally define $K_{i,j}\subset G$ by
    $\overline{x_{i,j}}K_{i,j}=\overline{x_{i,j}}K_i\setminus
    A_{i,j-1}$, and note by~\eqref{eq:weiss:1} that
    \begin{equation}\label{eq:weiss:9}
      \#(\overline{x_{i,j}}K_{i,j})\ge(1-\delta)\#K_i.
    \end{equation}

    If $\alpha_{i-1}\ge1$, then for all $j<i-1$ set
    $A_{j,0}=A_{i-1,0}$, $s(j)=0$ and $\nu_j=\nu_{i-1}$, and stop;
    otherwise, decrease $i$ and continue.

  After all these steps, we have obtained a decomposition
  \begin{equation}\label{eq:weiss:4}
    \Omega=\overline Q\sqcup\bigsqcup_{i=1}^t\bigsqcup_{j=1}^{s(i)}\overline{x_{i,j}}K_{i,j},
  \end{equation}
  with $\overline Q=\Omega\setminus A_{0,0}$. If the iteration was
  stopped because $\alpha_{i-1}\ge1$, then
  $\nu_0=\alpha_{i-1}(1-\delta)/\zeta>1-\epsilon/(2\#K)$; otherwise,
  $\nu_0\ge\overline{\nu_t}>1-\epsilon/(2\#K)$. In all cases, we have
  \begin{equation}\label{eq:weiss:7}
    \#\overline Q\le(1-\nu_0)\#\Omega<\frac\epsilon{2\#K}\#\Omega.
  \end{equation}
  Lift $\overline Q$ and $\overline{x_{i,j}}$ to $Q\subseteq G$ and
  $x_{i,j}\in G$ respectively. We have obtained a finite subset
  \begin{equation}\label{eq:weiss:5}
    T_n=Q\sqcup\bigsqcup_{i=1}^t\bigsqcup_{j=1}^{s(i)}x_{i,j}K_{i,j}
  \end{equation}
  of $G$.  Furthermore, the natural restriction
  $\pi\downharpoonright_{T_n}:T_n\to\Omega$ is a bijection, i.e.\
  $T_n$ is a transversal to $N_n$.  We compute
  \begin{align*}
    \#(T_n\cup T_nK)&=\#(T_nK)\text{ since }1\in K\\
    &\le\#(QK)+\sum_{i,j}\#(x_{i,j}K_iK)\text{ by }\eqref{eq:weiss:5}\\
    &\le\#Q\#K+\sum_{i,j}(1+{\textstyle\frac\epsilon2})(1-\delta)\#(x_{i,j}K_i)\text{ by }\eqref{eq:weiss:6}\\
    &\le\#Q\#K+\sum_{i,j}(1+{\textstyle\frac\epsilon2})\#(x_{i,j}K_{i,j})\text{ by }\eqref{eq:weiss:9}\\
    &\le\frac{\epsilon\#\Omega}{2\#K}\#K+(1+{\textstyle\frac\epsilon2})\#\Omega\text{ by }(\ref{eq:weiss:7}+\ref{eq:weiss:4})\\
    &\le(1+\epsilon)\#T_n.\qedhere
  \end{align*}
\end{proof}

\section{Tileable amenable algebras}\label{ss:tae}
We prove in this section an analogue of Theorem~\ref{thm:weiss} for
algebras. We follow as much as possible the notation of the previous
section; this mainly amounts to replacing `$\#$' by `$\dim$',
`$\sqcup$' by `$\oplus$' and so on. The equation numbers match between
these two sections; Section~\ref{ss:tag} also contains a few informal
remarks that may help the reader along the proof(s).

Throughout this section we consider an associative algebra $R$ over a
field $\K$, we use $\oplus$ and $\otimes$ to denote direct sum and
tensor products as $\K$-vector spaces, and we denote by `$\dim$' the
dimension as a $\K$-vector space. An \emph{i-subspace} of $R$ is a
subspace of $R$ admitting a basis consisting of invertible elements.
An \emph{i-algebra} is an algebra $R$ which is an i-subspace of
itself, and $R$ is \emph{i-amenable} if for every finite-dimensional
subspace $K$ and every $\epsilon>0$ there exists a
$(K,\epsilon)$-invariant i-subspace of $R$, in the sense of an
i-subspace $F$ such that
\[\frac{\rank(F+FK)-\rank(F)}{\rank(F)}<\epsilon.\]

\begin{lem}
  The group ring of an amenable group is an i-amenable i-algebra.
\end{lem}
\begin{proof}
  The basis $G$ of $\K G$ consists of invertible elements. Following
  the proof of Proposition~\ref{prop:crystal}, there exist
  $(K,\epsilon)$-invariant subspaces of the form $\K F$ with
  $F\subseteq G$, which have a basis $F$ consisting of invertible
  elements.
\end{proof}

Given a finite-dimensional subspace $A$ and a finite-dimensional
i-subspace $K$ of $R$, with fixed basis $\{k_1,\dots,k_d\}\subseteq
R^\times$, define
\[AK^*:=A\,\K\{k_1^{-1},\dots,k_d^{-1}\}.\]
Note that, although the notation does not make it apparent, $AK^*$ may
depend on the choice of basis of $K$.

\begin{thm}\label{thm:tiling}
  Let $R$ be an amenable i-algebra, let $K\le R$ be a
  finite-dimensional i-subspace, let $\epsilon>0$ be given, and let
  $I_0,I_1,\dots$ be a sequence of finite-codimension ideals in $R$
  such that $\bigcap_{n\in\N} I_n=\{0\}$.

  Then for all $n\gg0$ there exists a $(K,\epsilon)$-invariant subspace
  $T_n\le R$ that is a vector-space complement for $I_n$ in $R$.
\end{thm}

Let $\zeta\ge1$ and $\delta\in(0,1)$ be constants to be fixed later.
In this section we shall use, for $\mu\ge\delta$, the transformation
$\Theta_\mu:\R_+^2\to\R_+^2$ given by
\[\Theta_\mu(\nu,\alpha)=\Big(\nu+\mu(1-\alpha),\alpha+\frac{\mu(1-\alpha)}{1-\delta}\zeta\Big).\]
The proof of Theorem~\ref{thm:tiling} relies on the following
\begin{lem}\label{lem:tiling}
  Let $\Omega$ be a finite-dimensional $R$-module over an i-algebra
  $R$, with spanning subset $\mathscr C$; let $B$ be a subspace of
  $\Omega$ and let $K,L$ be finite-dimensional i-subspaces of $R$;
  assume that for all $x\in\mathscr C$ the orbit map
  \begin{equation}\label{eq:tiling:i}
    K\to\Omega,\; k\mapsto xk\qquad\text{is injective}.
  \end{equation}
  Let
  $\alpha,\nu\in[0,1)$ and $\zeta\ge1$ be such that
  \begin{xalignat*}{2}
    \dim(KL^*)&\le\zeta\dim K, & \dim(BK^*)&\le\alpha\dim\Omega,\\
    \dim B&=\nu\dim\Omega, & \dim(BL^*)&\le\alpha\dim\Omega;
  \end{xalignat*}
  Let furthermore $\delta\in(0,1)$ be given.

  Then there exist $s\ge1$, elements $x_1,\dots,x_s\in\mathscr C$, and
  $\mu\ge\delta$ such that, setting $B_0=B$ and $B_i=B_{i-1}+x_iK$ for
  $1\le i\le s$ and $\Theta_\mu(\nu,\alpha)=(\nu',\alpha')$, we have
  \begin{gather}
    \label{eq:tiling:1}
    \dim(x_iK\cap B_{i-1}) \le\delta\dim K\text{ for all }i\in\{1,\dots,s\};\\
    \label{eq:tiling:2} \dim B_s=\nu'\dim\Omega;\\
    \label{eq:tiling:3} \dim(B_sL^*)\le\alpha'\dim\Omega.
  \end{gather}
\end{lem}
\begin{proof}
  Let $x_1,\dots,x_s\in\mathscr C$ be a maximal-length sequence of
  elements such that~\eqref{eq:tiling:1} holds.  Let
  $\{v_1,\dots,v_d\}$ be a maximal subset of $\mathscr C$ whose image
  is independent in $\Omega/BK^*$; then $V:=\Bbbk\{v_1,\dots,v_d\}$ is
  a vector space complement to $BK^*$ in $\Omega$. By maximality of
  $(x_1,\dots,x_s)$, for all $x\in\mathscr C$ we have
  \begin{equation}\label{eq:tiling:p1}
    \dim(B_s\cap xK)=\dim\K\big\{x\otimes k\in\Omega\otimes K:\,xk\in B_s\big\}>\delta\dim K.
  \end{equation}
  
  Let $K$ have i-basis $\{k_i\}\subseteq R^\times$; then for any
  $A\le\Omega$ the vector space $\{\sum\omega_i\otimes k_i\in\Omega\otimes
  K:\,\omega_ik_i\in A\text{ for all }i\}$ is isomorphic to $A\otimes K$
  via the map
  \begin{equation}\label{eq:tensoriso}\tag{\ref{ss:tae}.\hbox{$\dagger$}}
    \sum \omega_i\otimes k_i\mapsto\sum \omega_ik_i\otimes k_i
  \end{equation}
  with inverse $\sum\xi_i\otimes k_i\mapsto\sum\xi_ik_i^{-1}\otimes
  k_i$.  We compute
  \begin{align*}
    \dim B_s\dim K&=\dim(B_s\otimes K)\\
    &=\dim\Big\{\sum \omega_i\otimes k_i\in\Omega\otimes
    K:\,\omega_i\in B_s\;\forall i\Big\}\\
    &=\dim\Big\{\sum\omega_i\otimes k_i\in\Omega\otimes K:\,\omega_ik_i\in
    B_s\;\forall i\Big\}\text{ by }\eqref{eq:tensoriso}\\
    &\ge\dim\Big\{\sum\omega_i\otimes k_i\in BK^*\otimes K:\,\omega_ik_i\in
    B_s\;\forall i\Big\}\\
    \intertext{\hfill $\displaystyle+\dim\Big\{\sum \omega_i\otimes
      k_i\in V\otimes K:\,\omega_ik_i\in B_s\;\forall i\Big\}$} 
    &\ge\dim\Big\{\sum\omega_ik_i^{-1}\otimes k_i\in BK^*\otimes
    K:\,\omega_i\in B_s\;\forall i\Big\}+\dim(VK\cap B_s)\\
    &\ge\dim(B\otimes K)+\sum_{j=1}^d\dim(v_jK\cap B_s)\text{ by }\eqref{eq:tiling:p1}\\
    &\ge\dim B\dim K+d\delta\dim K;
  \end{align*}
  and $d\ge(1-\alpha)\dim\Omega$. If we divide by $\dim K\dim\Omega$
  and set, $\mu=(\dim B_s/\dim\Omega-\nu)/(1-\alpha)$, then we get
  $\mu\ge\delta$, and therefore $s\ge1$. Note that~\eqref{eq:tiling:2}
  holds by the choice of $\mu$. Next, 
  \begin{align*}
    \dim B_s&\ge\dim B+\sum_{i=1}^s\dim(x_iK/(x_iK\cap B_{i-1}))\\
    &\ge\dim B+s(1-\delta)\dim K\text{ by }(\ref{eq:tiling:i},\ref{eq:tiling:1}),
  \end{align*}
  so $s\dim K\le(\nu'-\nu)\dim\Omega/(1-\delta)$; and
  \[\dim(B_sL^*)\le\dim(BL^*)+\sum_{i=1}^s\dim(x_iKL^*)
  \le\alpha\dim\Omega+s\zeta\dim K,
  \]
  from which~\eqref{eq:tiling:3} follows.
\end{proof}

For a sketch of the proof of Therom~\ref{thm:tiling}, which applies as
well to algebra setting as to the group setting, see
Remark~\ref{rem:tiling} on page~\pageref{rem:tiling}.

\begin{proof}[Proof of Theorem~\ref{thm:tiling}]
  Without loss of generality (possibly adjoining a unit to $R$ first),
  we suppose $1\in K$. Fix in order the following data:
  \begin{enumerate}
  \item Choose $\delta>0$, the ``allowed overlap'' as
    in~\eqref{eq:tiling:1}, such that $\delta\dim K<\frac\epsilon2$ and
    $(1+\frac\epsilon2)(1-\delta)>1$;
  \item Choose $\zeta>1$ such that $(1-\delta)/\zeta>1-\epsilon/(2\dim
    K)$, the ``relative F\o lner constant''.
  \end{enumerate}
  For all $t\in\N$ set
  $(\overline{\nu_t},\overline{\alpha_t})=\Theta_\delta^t(0,0)$.  Then
  it is clear that
  $\overline{\nu_t}/\overline{\alpha_t}=(1-\delta)/\zeta$ for all
  $t>0$, and that $\limsup\overline{\alpha_t}\ge$, so
  $\limsup\overline{\nu_t}\ge(1-\delta)/\zeta$.
  \begin{enumerate}\setcounter{enumi}{2}
  \item Let $t\in\N$, the ``tower height'', be such that
    $\overline{\nu_t}>1-\epsilon/(2\dim K)$.
  \end{enumerate}
  Then for any $\mu_1,\dots,\mu_t\ge\delta$ we also have
  \[\Theta_{\mu_t}(\cdots\Theta_{\mu_1}(0,0)\cdots)\in(1-\epsilon/(2\dim K),1]\times[0,1].\]
  \begin{enumerate}\setcounter{enumi}{3}
  \item Using i-amenability of $R$, construct finite-dimensional
    i-subspaces $K_0=K$, $K_1,\dots,K_t$ of $R$, the ``Rokhlin tower'',
    such that for all $j<i$ in $\{1,\dots,t\}$ one has
    $\dim(K_iK_j^*)<\zeta\dim K_i$, and such that
    \begin{equation}\label{eq:tiling:6}
      \dim(K_iK)\le\big(1+\frac\epsilon2\big)(1-\delta)\dim K_i \text{ for all
      }i\in\{1,\dots,t\}.
    \end{equation}
  \item Choose $n\in\N$, the ``quotient index'', large enough so that
    $K_iK_i^*\cap I_n=(0)$ for all $i\in\{1,\dots,t\}$.
  \end{enumerate}
  Set $\Omega=G/I_n$, let $\pi:R\to\Omega$ denote the natural quotient
  map, and let $\mathscr C$ denote the image in $\Omega$ of the
  spanning set $R^\times$ of $R$. Remark, by the choice of $n$, that
  $K_i\to\Omega$, $k\mapsto xk$ is injective for all
  $i\in\{1,\dots,t\}$ and all $x\in\mathscr C$.

  Now start with $A_{t,0}=\{0\}\le\Omega$ and
  $(\nu_t,\alpha_t)=(0,0)=(\overline{\nu_0},\overline{\alpha_0})$, and
  apply $t$ times Lemma~\ref{lem:tiling}: at step $i=t,t-1,\dots,1$,
  apply it with $B:=A_{i,0}$, $K:=\pi(K_i)$, $L:=\pi(K_{i-1})$, and
  $(\nu,\alpha):=(\nu_i,\alpha_i)$.  Lemma~\ref{lem:tiling} produces
  for us
  \begin{itemize}
  \item constants $s\ge1$ and $\mu\ge\delta$, which we rechristen
    $s(i)$ and $\mu_i$;
  \item a sequence $x_1,\dots,x_{s(i)}$ in $\mathscr C$, which we
    rechristen $\overline{x_{i,1}},\dots,\overline{x_{i,s(i)}}$;
  \item a sequence $B_0\subseteq\dots\subseteq B_{s(i)}$ of subspaces
    of $\Omega$, which we rechristen
    $A_{i,0},\dots,A_{i,s(i)}=:A_{i-1,0}$.
  \end{itemize}
  Set $(\nu_{i-1},\alpha_{i-1})=\Theta_{\mu_i}(\nu_i,\alpha_i)$, and
  note that $\nu_{i-1}\ge\overline{\nu_{t-(i-1)}}$ by induction on $i$
  and because $\mu_i\ge\delta$ and $\Theta_\mu$ is monotonically
  increasing in $\mu$.  Finally let $K_{i,j}\le K_i$ be any subspace
  such that $\overline{x_{i,j}}K_{i,j}\oplus(\overline{x_{i,j}}K_i\cap
  A_{i,j-1})=\overline{x_{i,j}}K_i$, and note by~\eqref{eq:tiling:1}
  that
  \begin{equation}\label{eq:tiling:9}
    \dim(\overline{x_{i,j}}K_{i,j})\ge(1-\delta)\dim K_i,
  \end{equation}
  because the $\overline{x_{i,j}}$ are invertible.

  If $\alpha_{i-1}\ge1$, then for all $j<i-1$ set $A_{j,0}=A_{i-1,0}$,
  $s(j)=0$ and $\nu_j=\nu_{i-1}$, and stop; otherwise, decrease $i$
  and continue.

  After all these steps, we have obtained a decomposition
  \begin{equation}\label{eq:tiling:4}
    \Omega=\overline Q\oplus\bigoplus_{i=1}^t\bigoplus_{j=1}^{s(i)}
    \overline{x_{i,j}}K_{i,j},
  \end{equation}
  where $\overline Q$ is any vector space complement to $A_{0,0}$ in
  $\Omega$. If the iteration was stopped because $\alpha_{i-1}\ge1$,
  then $\nu_0=\alpha_{i-1}(1-\delta)/\zeta>1-\epsilon/(2\dim K)$;
  otherwise, $\nu_0\ge\overline{\nu_t}>1-\epsilon/(2\dim K)$. In all
  cases, we have
  \begin{equation}\label{eq:tiling:7}
    \dim\overline Q\le(1-\nu_0)\dim\Omega
    <\frac\epsilon{2\dim K}\dim\Omega.
  \end{equation}
  Lift $\overline Q$ and $\overline{x_{i,j}}$ to $Q\le R$ and
  $x_{i,j}\in R^\times$ respectively.  We have obtained a
  finite-dimensional subspace
  \begin{equation}\label{eq:tiling:5}
    T_n=Q\oplus\bigoplus_{i=1}^t\bigoplus_{j=1}^{s(i)} x_{i,j}K_{i,j}
  \end{equation}
  of $R$. Furthermore, the natural restriction
  $\pi\downharpoonright_{T_n}:T_n\to\Omega$ is a bijection, i.e.\
  $T_n$ is a vector space complement to $I_n$. We compute
  \begin{align*}
    \dim(T_n+T_nK)&=\dim(TK)\text{ because }1\in K\\
    &\le\dim(QK)+\sum_{i,j}\dim(x_{i,j}K_iK)\text{ by }\eqref{eq:tiling:5}\\
    &\le\dim Q\dim
    K+\sum_{i,j}(1+\textstyle{\frac\epsilon2})(1-\delta)\dim(x_{i,j}K_i)\text{
      by }\eqref{eq:tiling:6}\\
    &\le\dim Q\dim K+\sum_{i,j}(1+\textstyle{\frac\epsilon2})\dim(x_{i,j}K_i)\text{ by }\eqref{eq:tiling:9}\\
    &\le\frac{\epsilon\dim\Omega}{2\dim K}\dim
    K+(1+\textstyle{\frac\epsilon2})\dim\Omega\text{ by }(\ref{eq:tiling:7}+\ref{eq:tiling:4})\\
    &\le(1+\epsilon)\dim T_n.\qedhere
  \end{align*}
\end{proof}

\section{Proof of Theorem~\ref{thm:vershik} and Corollary~\ref{cor:lcs}}\label{ss:vershik}
We start by a ``Reidemeister-Schreier'' result for algebras and ideals:
\begin{lem}\label{lem:rs}
  Let $R$ be a unital $\Bbbk$-algebra generated by a subspace $S$; let
  $I\triangleleft R$ be a right ideal in $R$; and let $F\le R$ be a
  complement of $I$, so we have $R=I\oplus_\Bbbk F$.  Let
  $x\mapsto\overline x$ be the projection $R\to F$; assume that $1\in
  F$ and $\overline 1=1$. Then
  \[I=\{fs-\overline{fs}:\, f\in F,\,s\in S\}R.\]
\end{lem}
\begin{proof}
  Write $J=\langle fs-\overline{fs}\rangle$; then obviously
  $J\subseteq I$. Conversely, consider first $s_1\dots s_n\in R$, and write
  $\equiv_J$ for congruence modulo $J$. Then
  \begin{align*}
    s_1\dots s_n &= (\overline1s_1-\overline{s_1})s_2\dots s_n+\overline{s_1}s_2\dots s_n\\
    &\equiv_J (\overline{s_1}s_2-\overline{\overline{s_1}s_2})s_3\dots s_n+\overline{\overline{s_1}s_2}s_3\dots s_n\\
    &\equiv_J\dots\equiv_J\overline{\overline{\overline{s_1}\cdots} s_n}.
  \end{align*}
  Consider now any $x=\sum_i s_{i,1}\dots s_{i,n_i}\in R$. Then
  $x\equiv_J\sum_i\overline{\overline{\overline{s_{i,1}}\cdots}s_{i,n_i}}\in
  F$, so $R=J+F$ and therefore $I=J$.
\end{proof}

\begin{rem}
  Consider the right ideal $I=(H-1)\Bbbk G\triangleleft\Bbbk G$ for
  some subgroup $H$ of $G=\langle S\rangle$. Let $T$ be a right
  transversal of $H$ in $G$; then $\Bbbk T$ is a complement of $I$ in
  $\Bbbk G$, so $I$ is generated by $\{fs-\overline{fs}:\,f\in T,s\in
  S\}$ by Lemma~\ref{lem:rs}, and therefore also by
  $\{fs\overline{fs}^{-1}-1\}$; so $H$ is generated by
  $\{fs\overline{fs}^{-1}:\,s\in S,f\in T\}$, which is the
  Reidemeister-Schreier generating set of $H$.\hfill $\triangle$
\end{rem}

\begin{cor}\label{cor:rs}
  Let $R=\langle S\rangle$ be an augmented algebra, and consider
  $I\triangleleft R$ with $R=I\oplus_\Bbbk F$. Then $I/I\varpi$ is
  spanned by the image of $(F+FS)\cap I$ in $I/I\varpi$.
\end{cor}

\noindent We are now ready to prove Theorem~\ref{thm:vershik}. We
start by a special case:
\begin{proof}[Proof of Theorem~\ref{thm:vershik} for $\Bbbk=\K$ a
  field of positive characteristic]
  Let a finite subset $S$ of $G$ and $\eta>1$ be given; we will show
  that $\lim\sqrt[n]{\dim(\K G/\varpi^n)}\le\eta$.

  Since $\K$ has positive characteristic and $G$ is finitely
  generated, the dimension subgroups $G_n=(1+\varpi^n)\cap G$ have
  finite index in $G$ for all $n$. Since any quotient of an amenable
  group is amenable, we may replace $G$ by $G/\bigcap_{n\ge1} G_n$ and
  assume from now on that $\bigcap G_n=1$. The ideals $\varpi^n$ have
  finite codimension in $\K G$, and have trivial intersection.

  We apply Theorem~\ref{thm:tiling} to $I_n=\varpi^n$: let $n_0\in\N$
  be such that for all $n\ge n_0$ there is a subspace $F_n$ of $\K G$
  with $\K G=I_n\oplus_\K F_n$ and $\dim(F_n+F_nS)<\eta\dim F_n$. Then
  \begin{align*}
    \dim(\K G/\varpi^{n+1})&=\dim(\K G/\varpi^n)+\dim(\varpi^n/\varpi^{n+1})\\
    &=\dim F_n+\dim(I_n/I_n\varpi)\\
    &\le\dim F_n+\dim((F_n+F_nS)\cap I_n)\text{ by Corollary~\ref{cor:rs}}\\
    &\le\dim((F_n+F_nS)\cap F_n)+\dim((F_n+F_nS)\cap I_n)\\
    &\le\dim((F_n+F_nS)\cap (F_n+I_n))=\dim(F_n+F_nS)\\
    &\le\eta\dim F_n=\eta\dim(\K G/\varpi^n).
  \end{align*}
  Set $C=\dim(\K G/\varpi^{n_0})/\eta^{n_0}$. We therefore have
  $\dim(\K G/\varpi^n)<C\eta^n$ for all $n\ge n_0$, so
  $\lim\sqrt[n]{\dim(\K G/\varpi^n)}\le\eta$ for all $\eta>1$.
\end{proof}

\begin{proof}[Proof of Theorem~\ref{thm:vershik} for $\Bbbk=\Z$]
  For $n\ge0$ set $r_n=\rank(\varpi^n/\varpi^{n+1})$, where $\varpi$
  denotes the augmentation ideal of $\Z G$; for $n\ge0$ and $p$ prime
  set $s_{n,p}=\rank(\varpi_p^n/\varpi_p^{n+1})$, where $\varpi_p$
  denotes the augmentation ideal of $\F G$; let
  $r_{n,p}=\rank(\varpi^n/\varpi^{n+1}\otimes_\Z\Z[\frac1p])$ denote
  the rank of the part of $\varpi^n/\varpi^{n+1}$ that is coprime to
  $p$, and let $r_{n,0}=\rank(\varpi^n/\varpi^{n+1}\otimes_\Z\Q)$
  denote the free rank of $\varpi^n/\varpi^{n+1}$.

  Since $G$ is finitely generated (say by $d$ elements),
  $\varpi^n/\varpi^{n+1}$ is a finite-rank abelian group (of rank at
  most $d^n$), so $\cong\Z^{r_{n,0}}\oplus\text{torsion}$. We thus
  have $r_{n,p}\ge r_{n,0}$, and for fixed $n$ we have $r_{n,p}=r_{n,0}$
  for almost all $p$.

  By Theorem~\ref{thm:vershik} for $\Bbbk=\F$, the sequence $s_{n,p}$
  grows subexponentially for any fixed $p$. Every $\F$ factor in
  $\varpi_p^n/\varpi_p^{n+1}$ lifts to a $\Z$-factor in $\Z G$, which
  then gives either a $\Z$-factor in $\varpi^n/\varpi^{n+1}$, or gives
  a torsion factor in $\varpi^m/\varpi^{m+1}$ for all $m\ge n$.
  Therefore $r_{m,p}\le\sum_{n\le m}s_{m,p}$, so $r_{n,p}$ grows
  subexponentially for any fixed $p$.

  The multiplication maps
  $\varpi^m/\varpi^{m+1}\otimes\varpi^n/\varpi^{n+1}\to\varpi^{m+n}/\varpi^{m+n+1}$
  are onto, so the sequence $r_n$ is submultiplicative ($r_mr_n\ge
  r_{m+n}$); the same holds for the sequences $r_{n,p}$ for fixed $p$,
  and for the sequence $r_{n,0}$.

  Let $\eta>1$ be given.  Then for some $n\in\N$ we have
  $r_{n,0}<\eta^n$; and for some $p_0$ we have $r_{n,p}=r_{n,0}$ if
  $p>p_0$. By submultiplicativity, $r_{kn,p}\le r_{n,p}^k<\eta^{kn}$
  for all $k$.

  For all $p\le p_0$ there exists $k_p\in\N$ such that
  $r_{k_pn,p}<\eta^{k_pn}$, because $r_{n,p}$ grows
  subexponentially. Set
  $m=n\cdot\operatorname{lcm}(k_1,\dots,k_{p_0})$. Then
  $r_{m,p}\le\eta^m$ for all $p$.

  We have $r_n=\max_pr_{n,p}$ for all $n\in\N$, so
  $r_{km}\le\eta^{km}$ for all $k\in\N$. Since $r_n$ is
  submultiplicative, we have
  \[\limsup\sqrt[n]{r_n}=\lim(r_{kn})^{1/kn}\le\eta\]
  \setbox9=\hbox{\cite{fekete:verteilung}*{page~233}}%
  by Fekete's Lemma~\cite{polya-s:analysis}*{volume~1, part~I,
    problem~98; originally~\box9}.
  Since $\eta>1$ was arbitrary, the sequence $r_n$ grows
  subexponentially.
\end{proof}

\begin{proof}[Proof of Theorem~\ref{thm:vershik} for general $\Bbbk$]
  Let $\varpi$ denote the augmentation ideal of $\Z G$, and let
  $\overline\varpi$ denote the augmentation ideal of $\Bbbk G$. Since
  the natural map
  $\varpi^n/\varpi^{n+1}\otimes\Bbbk\to\overline\varpi^n/\overline\varpi^{n+1}$
  is onto for all $n\in\N$, we have $\rank_\Bbbk(\Bbbk
  G/\overline\varpi^n)\le\rank_\Z(\Z G/\varpi^n)$, so the claim
  follows from Theorem~\ref{thm:vershik} for $\Bbbk=\Z$.
\end{proof}

\begin{rem}\label{rem:vershik}
  The converse of Theorem~\ref{thm:vershik} does not hold: the group
  $SL(d,\Z)$ for $d\ge3$ is certainly not amenable (it contains free
  subgroups), and neither is its congruence subgroup
  $K=\ker(SL(d,\Z)\to SL(d,\Z/p\Z))$. This subgroup is residually-$p$
  if $p\ge3$, since\footnote{As usual in these situations, one should
    treat sometimes $2$, sometimes $4$ as the even prime to extend
    this result to characteristic $2$.}  the subgroups
  $K_n=\ker(SL(d,\Z)\to SL(d,\Z/p^n\Z))$ have trivial intersection and
  index $p^{(n-1)(d^2-1)}$ in $K$. Then, because $K$ has the
  congruence property~\cite{bass-l-s:congruence}, the congruence
  subgroup $K_n$ coincides with the dimension subgroup as defined
  in~\S\ref{ss:augmented}, so $\dim(\varpi^n/\varpi^{n+1})$ grows
  subexponentially (approximately at rate
  $e^{(d^2-1)\pi\sqrt{2n/3}}$).  I wish to thank M.\ Ershov for
  pointing out this example to me.\hfill $\triangle$
\end{rem}

Corollary~\ref{cor:lcs} could follow along the same lines as the proof
for $\Bbbk=\Z$ of Theorem~\ref{thm:vershik}, by reducing from
$\gamma_n(G)/\gamma_{n+1}(G)$ to quotients of $p$-dimension subgroups
$G_{n,p}/G_{n+1,p}$ at all primes $p$, and using the fact that
$\bigoplus_{n\ge1}G_{n,p}/G_{n+1,p}$ is the primitive part of the Hopf
algebra $\overline{\F G}$ and therefore has subexponential growth. We
will however opt for a shortcut:
\begin{proof}[Proof of Corollary~\ref{cor:lcs}]
  The classical \emph{dimension subgroups} of $G$ are the subgroups
  $\delta_n(G)=G\cap(1+\varpi^n)$, where $\varpi$ denotes the
  augmentation ideal in $\Z G$. By a result of
  Gupta~\cite{gupta:odddimension}, the quotient
  $\delta_n(G)/\gamma_n(G)$ is a finite $2$-group. Now
  $\rank(\delta_n(G)/\delta_{n+1}(G))$ grows subexponentially by
  Theorem~\ref{thm:vershik} for $\Bbbk=\Z$, since
  $\delta_n(G)/\delta_{n+1}(G)$ is a submodule of
  $\varpi^n/\varpi^{n+1}$; and $\dim_{\F[2]}G_{n,2}/G_{n+1,2}$ grows
  subexponentially since $\overline{\F[2]G}$ has subexponential growth
  by Theorem for $\Bbbk=\F[2]$. We conclude that
  \begin{align*}
    \rank(\gamma_n(G)/\gamma_{n+1}(G))&\le\rank(\gamma_n(G)/\gamma_{n+1}(G)\otimes\Z[{\textstyle\frac12}])+\rank(\gamma_n(G)/\gamma_{n+1}(G)\otimes\F[2])\\
    &\le\rank(\delta_n(G)/\delta_{n+1}(G))+\rank(G_{n,2}/G_{n+1,2})
  \end{align*}
  grows subexponentially.
\end{proof}

%


\bib{bartholdi-g:lie}{incollection}{
  author={Bartholdi, Laurent},
  author={Grigorchuk, Rostislav I.},
  title={Lie methods in growth of groups and groups of finite width},
  conference={ title={Computational and geometric aspects of modern algebra (Edinburgh, 1998)}, },
  book={ series={London Math. Soc. Lecture Note Ser.}, volume={275}, publisher={Cambridge Univ. Press}, place={Cambridge}, },
  date={2000},
  pages={1--27},
  review={\MR {1776763 (2001h:20046)}},
  eprint={arXiv.org/abs/math/0002010},
}

\bib{bartholdi:aa1}{article}{
  author={Bartholdi, Laurent},
  title={On amenability of group algebras, I},
  date={2007},
  status={to appear in Isr. J. Math.},
  eprint={arXiv.org/abs/math/0608302},
}

\bib{bass-l-s:congruence}{article}{
  author={Bass, Hyman},
  author={Lazard, Michel},
  author={Serre, Jean-Pierre},
  title={Sous-groupes d'indice fini dans ${\bf SL}(n,\,{\bf Z})$},
  language={French},
  journal={Bull. Amer. Math. Soc.},
  volume={70},
  date={1964},
  pages={385--392},
  review={\MR {0161913 (28 \#5117)}},
}

\bib{bogopolski:bilipschitz}{article}{
  author={Bogopol{\cprime }ski{\u \i }, Oleg~V.},
  title={Infinite commensurable hyperbolic groups are bi-Lipschitz equivalent},
  language={Russian, with Russian summary},
  journal={Algebra i Logika},
  volume={36},
  date={1997},
  number={3},
  pages={259--272, 357},
  issn={0373-9252},
  translation={ journal={Algebra and Logic}, volume={36}, date={1997}, number={3}, pages={155--163}, issn={0002-5232}, },
  review={\MR {1485595 (98h:57002)}},
}

\bib{burnside:question}{article}{
  author={Burnside, William},
  title={On an unsettled question in the theory of discontinuous groups},
  date={1902},
  journal={Quart. J. Pure Appl. Math.},
  volume={33},
  pages={230\ndash 238},
}

\bib{elek:amenaa}{article}{
  author={Elek, G{\'a}bor},
  title={The amenability of affine algebras},
  journal={J. Algebra},
  volume={264},
  date={2003},
  number={2},
  pages={469--478},
  issn={0021-8693},
  review={\MR {1981416 (2004d:16043)}},
}

\bib{ershov:goshat}{article}{
  author={Ershov, Mikhail V.},
  title={Golod-Shafarevich groups with property (T) and Kac-Moody groups},
  date={2006},
  note={Preprint},
}

\bib{fekete:verteilung}{article}{
  author={Fekete, Mih\'aly},
  title={\"Uber die Verteilung der Wurzeln bei gewissen algebraischen Gleichungen mit ganzzahligen Koeffizienten},
  language={German},
  journal={Math. Z.},
  volume={17},
  date={1923},
  number={1},
  pages={228--249},
  issn={0025-5874},
  review={\MR {1544613}},
}

\bib{folner:banach}{article}{
  author={F{\o }lner, Erling},
  title={Note on groups with and without full Banach mean value},
  date={1957},
  journal={Math. Scand.},
  volume={5},
  pages={5\ndash 11},
}

\bib{golod:nil}{article}{
  author={Golod, Evgueni{\u \i }~S.},
  title={On nil-algebras and finitely approximable $p$-groups},
  date={1964},
  journal={Izv. Akad. Nauk SSSR Ser. Mat.},
  volume={28},
  pages={273\ndash 276},
}

\bib{gromov:topinv1}{article}{
  author={Gromov, Mikhael L.},
  title={Topological invariants of dynamical systems and spaces of holomorphic maps. I},
  journal={Math. Phys. Anal. Geom.},
  volume={2},
  date={1999},
  number={4},
  pages={323--415},
  issn={1385-0172},
  review={\MR {1742309 (2001j:37037)}},
}

\bib{gupta:odddimension}{article}{
  author={Gupta, Narain~D.},
  title={The dimension subgroup conjecture holds for odd order groups},
  journal={J. Group Theory},
  volume={5},
  date={2002},
  number={4},
  pages={481--491},
  issn={1433-5883},
  review={\MR {1931371 (2003m:20019)}},
}

\bib{harpe:mfap}{article}{
  author={de la Harpe, Pierre},
  title={Mesures finiment additives et paradoxes},
  language={French, with English and French summaries},
  conference={ title={Autour du centenaire Lebesgue}, },
  book={ series={Panor. Synth\`eses}, volume={18}, publisher={Soc. Math. France}, place={Paris}, },
  date={2004},
  pages={39--61},
  review={\MR {2143416 (2006e:43001)}},
}

\bib{harpe:uniform}{article}{
  author={Harpe, Pierre~{de la}},
  title={Uniform growth in groups of exponential growth},
  date={2002},
  journal={Geom. Dedicata},
  volume={95},
  pages={1\ndash 17},
}

\bib{herstein:rings}{book}{
  author={Herstein, Israel~N.},
  title={Noncommutative rings},
  series={Carus Mathematical Monographs},
  publisher={Mathematical Association of America},
  address={Washington, DC},
  date={1994},
  volume={15},
  isbn={0-88385-015-X},
  note={Reprint of the 1968 original; With an afterword by Lance W. Small},
}

\bib{huppert-b:fg2}{book}{
  author={Huppert, Bertram},
  author={Blackburn, Norman},
  title={Finite groups. II},
  series={Grundlehren der Mathematischen Wissenschaften [Fundamental Principles of Mathematical Sciences]},
  volume={242},
  note={AMD, 44},
  publisher={Springer-Verlag},
  place={Berlin},
  date={1982},
  pages={xiii+531},
  isbn={3-540-10632-4},
  review={\MR {650245 (84i:20001a)}},
}

\bib{jennings:gpring}{article}{
  author={Jennings, Stephen~A.},
  title={The structure of the group ring of a $p$-group over a modular field},
  date={1941},
  journal={Trans. Amer. Math. Soc.},
  volume={50},
  pages={175\ndash 185},
}

\bib{koch:galois}{book}{
  author={Koch, Helmut},
  title={Galoissche Theorie der $p$-Erweiterungen},
  publisher={Springer-Verlag},
  address={Berlin},
  date={1970},
  note={Mit einem Geleitwort von I. R. \v Safarevi\v c},
}

\bib{lazard:uea}{article}{
  author={Lazard, Michel},
  title={Sur les alg\`ebres enveloppantes universelles de certaines alg\`ebres de Lie},
  date={1952},
  journal={C. R. Acad. Sci. Paris S{\'e}r. I Math.},
  volume={234},
  pages={788\ndash 791},
}

\bib{ornstein-weiss:entropyiso}{article}{
  author={Ornstein, Donald S.},
  author={Weiss, Benjamin},
  title={Entropy and isomorphism theorems for actions of amenable groups},
  journal={J. Analyse Math.},
  volume={48},
  date={1987},
  pages={1--141},
  issn={0021-7670},
  review={\MR {910005 (88j:28014)}},
}

\bib{petrogradsky:polynilpotent}{article}{
  author={Petrogradski{\u \i }, Victor~M.},
  title={Growth of finitely generated polynilpotent Lie algebras and groups, generalized partitions, and functions analytic in the unit circle},
  date={1999},
  issn={0218-1967},
  journal={Internat. J. Algebra Comput.},
  volume={9},
  number={2},
  pages={179\ndash 212},
  review={\MRhref {1 703 073}},
}

\bib{polya-s:analysis}{book}{
  author={Poly{\'a}, George},
  author={Szeg{\"o}, Gabor},
  title={Aufgaben und Lehrs\"atze aus der Analysis. Band I: Reihen. Integralrechnung. Funktionentheorie},
  language={German},
  series={Dritte berichtigte Auflage. Die Grundlehren der Mathematischen Wissenschaften, Band 19},
  publisher={Springer-Verlag},
  place={Berlin},
  date={1964},
  pages={xvi+338},
  review={\MR {0170985 (30 \#1219a)}},
}

\bib{vaillant:folner}{article}{
  author={Vaillant, Ghislain},
  title={Folner conditions, nuclearity, and subexponential growth in $C\sp *$-algebras},
  journal={J. Funct. Anal.},
  volume={141},
  date={1996},
  number={2},
  pages={435--448},
  issn={0022-1236},
  review={\MR {1418514 (97m:46098)}},
}

\bib{vershik:amenability}{article}{
  author={Vershik, Anatoly M.},
  title={Amenability and approximation of infinite groups},
  note={Selected translations},
  journal={Selecta Math. Soviet.},
  volume={2},
  date={1982},
  number={4},
  pages={311--330},
  issn={0272-9903},
  review={\MR {721030 (86g:43006)}},
}

\bib{weiss:monotileable}{article}{
  author={Weiss, Benjamin},
  title={Monotileable amenable groups},
  conference={ title={Topology, ergodic theory, real algebraic geometry}, },
  book={ series={Amer. Math. Soc. Transl. Ser. 2}, volume={202}, publisher={Amer. Math. Soc.}, place={Providence, RI}, },
  date={2001},
  pages={257--262},
  review={\MR {1819193 (2001m:22014)}},
}



\begin{bibsection}
\begin{biblist}
\bibselect{bartholdi,math}
\end{biblist}
\end{bibsection}
\end{document}